\documentclass[11pt]{article}

\usepackage[mathscr]{eucal}
\usepackage{color}
\usepackage{amsmath}
\usepackage{amssymb}
\usepackage{amsfonts}
\usepackage{amscd}
\usepackage{amsthm}
\usepackage{latexsym}
\usepackage{graphicx}
\usepackage{subfig}

\def\be{\begin{equation}}
\def\ee{\end{equation}}
\def\bse{\begin{subequations}}
\def\ese{\end{subequations}}

\newtheorem{thm}{Theorem}

\newtheorem{prop}[thm]{Proposition}

\usepackage{latexsym}

\let\er\eqref

\let\be\beta

\newcommand{\R}{{\mathbb R}}

\newtheorem{theorem}{Theorem}
\newtheorem{lemma}[theorem]{Lemma}

\def\bse{\begin{subequations}}
\def\ese{\end{subequations}}

\usepackage{geometry}

\title{A nonlocal reaction diffusion equation and its relation with Fujita exponent}

\author{ Shen Bian\footnote{Beijing University of Chemical Technology, 100029, Beijing. Email: \texttt{bianshen66@163.com}. Partially supported by National Science Foundation of China (Grant No. 11501025), China Postdoctoral Science Foundation (Grant No. 2014M560037) and the Fundamental Research Funds for the Central Universities (Grant No. ZY1528).}
  \and Li Chen\footnote{Universit\"at Mannheim, 68131, Mannheim. Email: \texttt{chen@math.uni-mannheim.de}. Partially supported by the National
Natural Science Foundation of China (NSFC), No. 11271218.} }
\date{}

\begin{document}
\let\cleardoublepage\clearpage

\maketitle

\begin{abstract}
This paper is concerned with a type of nonlinear reaction-diffusion equation, which arises from the population dynamics. The equation includes a certain type reaction term $u^\alpha(1- \sigma \int_{\R^n}u^\beta dx)$ of dimension $n \ge 1$ and $\sigma>0$. An energy-methods-based proof on the existence of global solutions is presented and the qualitative behavior of solution which is decided by the choice of $\alpha,\beta$ is exhibited. More precisely, for $1 \le \alpha<1+(1-2/p)\beta$, where $p$ is the exponent appears in Sobolev's embedding theorem defined in \er{p}, the equation admits a unique global solution for any nonnegative initial data. Especially, in the case of $n\geq 2$ and $\beta=1$, the exponent $\alpha<1+2/n$ is exactly the well-known Fujita exponent. The global existence result obtained in this paper shows that by switching on the nonlocal effect, i.e., from $\sigma=0$ to $\sigma>0$, the solution's behavior differs distinctly, that's, from finite time blow-up to global existence.
\end{abstract}

\section{Introduction}

In this paper, we study the following nonlocal initial boundary value problem,
\bse\label{nkpp}
\begin{align}
&\ u_t-\Delta u=u^\alpha \left(1- \sigma \int_{\R^n} u^\beta(x,t)dx\right),\quad   && x \in \R^n, t>0,
\\
&\ u(x,0)=u_0(x)\geq0,\quad &&x \in \R^n,
\end{align}
\ese where $u$ is the density, $\alpha,\beta \geq 1$, $\sigma>0$.

This model is developed to describe the population dynamics \cite{f37,kpp} with the form $\frac{\partial u}{\partial t}=\frac{\partial^2 u}{\partial x^2}+F(u)$, the function $F$ is considered as the rate of the reproduction. In this paper, we will study problems with nonlocal reaction term.

As appeared in many literatures \cite{BCL15,hy95,lcl09,Lorz:2011hl,Lorz:2013vp,vol1,vol2,VVpp,ww11,ww96,Yi96}, nonlocal type reaction terms can describe also Darwinian evolution of a structured population density or the behaviors of cancer cells with therapy. We review some of the known results on the reaction-diffusion equation with a nonlocal term. In \cite{bb82}, the authors considered the equation with reaction term
$F(t,u)=e^u+ \int_{\Omega} e^u dx$,
where $\Omega$ is a bounded domain of $\R^n$, for which the above problem represents an ignition model for a compressible reactive gas, and they proved the finite time blow-up of solutions.
Later on, a power-like nonlinearity $F(t,u)=\int_{\Omega} u^p(t,y) dy- k u^q(t,x)$ was investigated by Wang and Wang \cite{ww96} with $p,q>1$, and they proved that solutions blow up.
Moreover, \cite{hy95} studied the case $F(u)=u^p-\frac{1}{|\Omega|} \int_{\Omega} u^p(t,y)dy$, this typical structure has mass conservation, and the authors showed that if $p>n/(n-2)$, the solutions will blow up in finite time with any initial data, while for $1<p<n/(n-2)$, the solution exists globally.
Recently, the authors in \cite{BCL15} studied the initial boundary value problem of (\ref{nkpp}) for $\beta=1$ and $\sigma=1$ in a bounded domain $\Omega\subset\R^n$. Global existence, uniqueness and long time behavior of solution were obtained, where the whole estimates rely on the long time asymptotics of the total mass $\int_{\Omega} u(t,x) dx$. However this property is not valid in the whole space case. That gives us a motivation to study the whole space case.

Another reason to study the Cauchy problem for this nonlocal reaction diffusion equation is that it has a close relation to the well-known Fujita exponent. We first list the main result of this paper, afterwards explain its relation to the Fujita exponent.

We take $\sigma=1$ for simplicity. All the discussions and results that obtained in this paper are valid for any positive $\sigma$.

\begin{thm}\label{thm1}
Assume that $u_0$ is nonnegative and $u_0 \in L^\beta (\R^n)\cap L^\infty(\R^n)$, $n \geq 1$. If $\alpha$ satisfies
   \begin{align}
    1 \le \alpha<1+(1-2/p)\beta,
   \end{align}
where $p$ is the exponent from the Sobolev embedding theorem, i.e.
\begin{align}\label{p}
\left\{
    \begin{array}{ll}
      p=\frac{2n}{n-2}, & n \ge 3, \\
      2<p<\infty , & n=2, \\
      p=\infty, & n=1.
    \end{array}
  \right.
\end{align}
then problem (\ref{nkpp}) has a unique bounded nonnegative solution. Moreover, the following
a priori estimates hold true. That's for any $t>0$ and $\beta \le k \le \infty$
\begin{align}
\int_{\R^n}  u(t)^k dx \le C\left(\|u_0\|_{L^\infty(\R^n)},\|u_0\|_{L^\beta(\R^n)}\right).
\end{align}
\end{thm}

Fujita in \cite{fuj1} showed that in the case of $\sigma=0$, the Cauchy problem has no global solution for super-critical exponent $\alpha<1+2/n$ with any nonnegative initial data by using comparison principle. On the other hand, in case of $\alpha>1+2/n$, there exists a global solution for sufficiently small initial data, and no global solutions for sufficiently large initial data. Later on, the authors in \cite{kh73} considered the case $n=2,\alpha=2$ and proved that it has no global solution for any nontrivial nonnegative initial data.

Compared to the case $\sigma=0$, the main mathematical difficulty in studying our problem is that unlike the cases in \cite{fuj1,hy95,ww96}, solutions of \er{nkpp} with positive $\sigma$ do not obey the comparison principle and mass conservation, which makes the use of many technical conditions or tools impossible. In addition, if $1- \sigma \int_{\R^n} u^\beta dx$ remains positive, it has similar structure to the case $\sigma=0$, therefore our problem might have no global solution for $\alpha<1+2/n$. However, the results obtained in theorem \ref{thm1} give an opposite consequence. Furthermore, in case $\beta=1$, within exactly the same range for $\alpha$, i.e., $1<\alpha<1+2/n$, $\sigma=0$ gives the finite time blow up of the solution, while $\sigma>0$ gives always global existence of the solution. In other words, switching on the nonlocal effect before the blow up time will prevent the solution's blow-up behavior.

Next we give a brief outline of the key estimates in order to get global existence. Actually, most of the {\it a priori} estimates are based on the following arguments. For any $k \ge 1$
\begin{align}
\frac{d}{dt} \int_{\R^n}  u^k dx + \frac{4(k-1)}{k} \int_{\R^n}
|\nabla u^{\frac{k}{2}} |^2 dx+ k \int_{\R^n} u^\beta dx \int_{\R^n} u^{k+\alpha-1}
dx = k \int_{\R^n} u^{k+\alpha-1} dx. \label{guji1}
\end{align}
Because of the speciality of the exponent $\beta+\alpha-1$, which means
\begin{align}
\|u\|_{L^{\beta+\alpha-1}(\R^n)}^{\beta+\alpha-1} & \le  \left( \|u\|_{L^{\beta+2(\alpha-1)}(\R^n)}^{\beta+2(\alpha-1)} \|u\|_{L^\beta(\R^n)}^{\beta} \right)^{\frac{1}{2}},
\end{align}
we can get the estimate for $L^{\beta+\alpha-1}$ norm from \eqref{guji1} in the first stage. Afterwards, the estimates for $L^k$ norm (with finite $k$) of solutions will be splitted into two cases: $\beta \le k \le \beta+\alpha-1$ and $k>\beta+\alpha-1$. More precisely, the uniform boundedness of $L^k$ norm of solutions for $\beta \le k \le \beta+\alpha-1$ is obtained by making use of the nonnegative term $\int_{\R^n} u^\beta dx \int_{\R^n} u^{k+\alpha-1}dx$ and the evolution of $L^\beta$ norm, where $L^\beta$ norm itself is an independent delicate case. While for large $\beta+\alpha-1<k<\infty$, $\int_{\R^n}
|\nabla u^{\frac{k}{2}} |^2 dx$ will be applied due to Sobolev's embedding theorem. In the end, the uniform in time $L^\infty$ norm is done by using a modified Moser type iteration argument.



\section{Global Existence for $1 \le \alpha < 1+(1-2/p)\beta$}
The main task in proving the global existence is to get the \emph{a priori} estimates in Proposition \ref{global}. Afterwards, a direct application of the standard parabolic theory leads to global existence of a unique solution. Before going to the proof, we need the following Sobolev inequality
\begin{lemma}[\cite{lieb202}] \label{sobolev}
For $n \ge 1$, $p$ is expressed by \er{p}, let $u \in H^1(\R^n)$. Then $u \in L^p(\R^n)$ and the following inequality holds:
\begin{align}
\|u\|_{L^p(\R^n)}^2 &\le C(n) \|\nabla u\|_{L^2(\R^n)}^2,~~n \ge 3, \\
\|u\|_{L^p(\R^n)}^2 &\le C(n) \left( \|\nabla u\|_{L^2(\R^n)}^2+\|u\|_{L^2(\R^n)}^2 \right),~~n=1,2.
\end{align}
Here $C(n)$ are given constants depending on $n$.
\end{lemma}
Based on the above Sobolev inequality, one has the following results which will be used in the proof of Proposition \ref{global}.
\begin{lemma}\label{interpolation}
Let $n \ge 1$. $p$ is expressed by \er{p}, $1 \le r<q<p$ and $\frac{q}{r}<\frac{2}{r}+1-\frac{2}{p}$, then for $v\in H^1(\R^n)$, it holds
\begin{align}
\|v\|_{L^q(\R^n)}^q \le C(n) C_0^{-\frac{ \lambda q}{2-\lambda q}} \|v\|_{L^r(\R^n)}^{\gamma} + C_0 \|\nabla v\|_{L^2(\R^n)}^2,~~n \ge 3, \label{nge3} \\
\|v\|_{L^q(\R^n)}^q \le C(n) \left(C_0^{-\frac{ \lambda q}{2-\lambda q}}  + C_1^{-\frac{ \lambda q}{2-\lambda q}} \right)\|v\|_{L^r(\R^n)}^{\gamma} + C_0 \|\nabla v\|_{L^2(\R^n)}^2+ C_1 \|v\|_{L^2(\R^n)}^2,~~n=1,2. \label{n12}
\end{align}
Here $C(n)$ are constants depending on $n$, $C_0,C_1$ are arbitrary positive constants and
\begin{align}
\lambda=\frac{\frac{1}{r}-\frac{1}{q}}{\frac{1}{r}-\frac{1}{p}} \in (0,1),~~
\gamma=\frac{2(1-\lambda) q}{2-\lambda q}=\frac{2\left(1-\frac{q}{p}\right)}{\frac{2-q}{r}-\frac{2}{p}+1}.
\end{align}
\end{lemma}
\noindent {\bf Proof.} First of all, it is easy to verify that $\lambda q<2$ if $\frac{q}{r}<\frac{2}{r}+1-\frac{2}{p}$. Therefore, by H\"{o}lder's inequality and Lemma \ref{sobolev} one has that for $n \ge 3$,
\begin{align*}
\|v\|_{L^q(\R^n)}^q & \le \|v\|_{L^p(\R^n)}^{\lambda q} \|v\|_{L^r(\R^n)}^{(1-\lambda)q} \\
& \le C(n) \|\nabla v\|_{L^2(\R^n)}^{\lambda q} \|v\|_{L^r(\R^n)}^{(1-\lambda)q}.
\end{align*}
Since $\lambda q<2$, then by Young's inequality one has
\begin{align*}
\|v\|_{L^q(\R^n)}^q \le C(n) C_0^{-\frac{ \lambda q}{2-\lambda q}} \|v\|_{L^r(\R^n)}^{\frac{2(1-\lambda) q}{2-\lambda q}} + C_0 \|\nabla v\|_{L^2(\R^n)}^2.
\end{align*}
For $n=1,2$, similar proof to the case $n \ge 3$, by H\"{o}lder's inequality and Young's inequality we obtain
\begin{align*}
\|v\|_{L^q(\R^n)}^q & \le \|v\|_{L^p(\R^n)}^{\lambda q} \|v\|_{L^r(\R^n)}^{(1-\lambda)q} \\
& \le C(n) \left( \|\nabla v\|_{L^2(\R^n)}^2+ \|v\|_{L^2(\R^n)}^2 \right)^{\frac{\lambda q}{2}} \|v\|_{L^r(\R^n)}^{(1-\lambda)q} \\
& \le C(n) \left( \|\nabla v\|_{L^2(\R^n)}^{\lambda q}+\|v\|_{L^2(\R^n)}^{\lambda q} \right) \|v\|_{L^r(\R^n)}^{(1-\lambda)q}\\
& \le C(n) \left(C_0^{-\frac{ \lambda q}{2-\lambda q}}  + C_1^{-\frac{ \lambda q}{2-\lambda q}} \right)\|v\|_{L^r(\R^n)}^{\frac{2(1-\lambda) q}{2-\lambda q}} + C_0 \|\nabla v\|_{L^2(\R^n)}^2+ C_1 \|v\|_{L^2(\R^n)}^2.
\end{align*}
Here, for simplicity, $C(n)$ are bounded constants depending on $n$.
This concludes the proof. $\Box$

Now we consider the a priori estimates for global existence.
\begin{prop}\label{global}
Let $n \geq 1$, $\alpha,\beta \ge 1$, $p$ is defined as in \er{p}. Assume $u_0$ is nonnegative and $u_0 \in L^\beta \cap L^\infty(\R^n)$. If $~\alpha$ satisfies
   \begin{align}
    1 \le \alpha<1+(1-2/p)\beta,
   \end{align}
then any nonnegative solution of (\ref{nkpp}) satisfies that for any $\beta \le k \le \infty$ and any $t>0$
\begin{align}
\int_{\R^n}  u(t)^k dx \le C(\alpha,\beta,K_0).
\end{align}
Here $K_0=\max\left\{ 1,\|u_0\|_{L^\infty(\R^n)},\|u_0\|_{L^\beta(\R^n)} \right\}$.
\end{prop}
\noindent\textbf{Proof of Proposition \ref{global}.} The proof will be given step by step. Firstly, we will give the a priori estimates for $L^k$ norm of solution for any $k>\max\left\{ \frac{2(\alpha-1)}{p-2},1 \right\}$. Then based on the a priori estimates, we will show that for any $\beta \le k \le \beta+\alpha-1$, the $L^k$ norm is uniformly bounded for any $t>0$. Thirdly, the boundedness of $L^k$ norm for $\beta+\alpha-1<k<\infty$ is proved. Finally, it follows that the $L^\infty$ norm of solutions is uniformly bounded by the iterative method. Thus closes the proof.

First if all, it is obtained by multiplying \er{nkpp} with $k u^{k-1}(k \ge 1)$
\begin{align}
\frac{d}{dt} \int_{\R^n}  u^k dx + \frac{4(k-1)}{k} \int_{\R^n}
|\nabla u^{\frac{k}{2}} |^2 dx+ k \int_{\R^n} u^\beta dx \int_{\R^n} u^{k+\alpha-1}
dx = k \int_{\R^n} u^{k+\alpha-1} dx. \label{guji2}
\end{align}
Noticing the two nonnegative terms of the left hand side of \er{guji2}, we firstly use the nonnegative term $\int_{\R^n} u^\beta dx \int_{\R^n} u^{k+\alpha-1} dx$ to control the right hand side of \er{guji2} in order to get the uniform boundedness of $L^k$ norm of solutions for some suitable small $k$. For $L^k$ norm of large $k$, we will take advantage of the other nonnegative term $\int_{\R^n} |\nabla u^{\frac{k}{2}} |^2 dx$ to dominate the right hand side of \er{guji2}. Finally, the iterative method is applied to prove the $L^\infty$ norm of solutions.

\noindent {\bf Step 1 (A priori estimates).} Taking $k>\max\left\{ \frac{2(\alpha-1)}{p-2},1 \right\}$ and $\max\left\{ \frac{(\alpha-1) p}{p-2},1 \right\}<k'<k+\alpha-1$ such that we can let $$
v=u^{\frac{k}{2}},~~q=\frac{2(k+\alpha-1)}{k},~~r=\frac{2k'}{k},~~C_0=\frac{k-1}{k^2}, ~~C_1=\frac{1}{2k}
$$
in Lemma \ref{interpolation}, this follows
\begin{align}\label{guji4}
\begin{array}{ll}
\int_{\R^n} u^{k+\alpha-1} dx &\le\frac{k-1}{k^2} \left\|\nabla
u^{\frac{k}{2}} \right\|_{L^2(\R^n)}^2 + C(k) \left\|
u^{\frac{k}{2}} \right\|_{L^\frac{2k'}{k}(\R^n)}^{\frac{2b}{k}},~~n \ge 3, \\
\int_{\R^n} u^{k+\alpha-1} dx &\le\frac{k-1}{k^2} \left\|\nabla
u^{\frac{k}{2}} \right\|_{L^2(\R^n)}^2 + C(k) \left\|
u^{\frac{k}{2}} \right\|_{L^\frac{2k'}{k}(\R^n)}^{\frac{2b}{k}}+\frac{1}{2k} \|u\|_{L^k(\R^n)}^k,~~n=1,2.
\end{array}
\end{align}
where
\begin{align}\label{b}
b=\frac{(1-\lambda)(k+\alpha-1)}{1-\frac{\lambda(k+\alpha-1)}{k}},~~\lambda=\frac{ \frac{k}{2k'}- \frac{k}{2(k+\alpha-1)} }{
\frac{k}{2k'} -\frac{1}{p} } \in (0,1).
\end{align}
Thus, (\ref{guji2}) with (\ref{guji4}) yields
\begin{align}
& \frac{d}{dt} \int_{\R^n}  u^k dx + k \int_{\R^n} u^\beta dx \int_{\R^n}
u^{k+\alpha-1} dx + \frac{3(k-1)}{k} \| \nabla u^{\frac{k}{2}} \|_{L^{2}(\R^n)}^2 \nonumber\\
\le ~ & C(k) \|u\|_{L^{k'}(\R^n)}^{b},~~n \ge 3. \label{nge3priori} \\
& \frac{d}{dt} \int_{\R^n}  u^k dx + k \int_{\R^n} u^\beta dx \int_{\R^n}
u^{k+\alpha-1} dx + \frac{3(k-1)}{k} \| \nabla u^{\frac{k}{2}} \|_{L^{2}(\R^n)}^2 \nonumber\\
\le ~ & C(k) \|u\|_{L^{k'}(\R^n)}^{b}+\frac{1}{2} \|u\|_{L^k(\R^n)}^k ,~~n=1,2, \label{n12priori}
\end{align}

Besides, with the help of interpolation inequality, in case of $\max\left\{ \frac{p(\alpha-1)}{p-2},\beta\right\}<k'<k+\alpha-1,$
\begin{align}
\|u\|_{L^{k'}(\R^n)}^{b} & \le \|u\|_{L^{k+\alpha-1}(\R^n)}^{b \theta} \|u\|_{L^\beta(\R^n)}^{(1-\theta)b} \nonumber\\
& = \left(  \|u\|_{L^{k+\alpha-1}(\R^n)}^{k+\alpha-1} \|u\|_{L^\beta(\R^n)}^\beta \right)^{\frac{b \theta}{k+\alpha-1}} \|u\|_{L^\beta(\R^n)}^{b\left(1-\theta-\frac{ \theta \beta}{k+\alpha-1} \right)}  \label{inter}
\end{align}
where
\begin{align}\label{theta}
\theta=\frac{\frac{1}{\beta}-\frac{1}{k'}}{\frac{1}{\beta}-\frac{1}{k+\alpha-1}} \in (0,1).
\end{align}
In addition, simple computations show that
\begin{align}\label{xiaoyu}
\frac{b \theta}{k+\alpha-1}<1
\end{align}
if and only if
\begin{align}\label{xiaoyup}
1\le \alpha<1+\left(1-\frac{2}{p}\right)\beta.
\end{align}
Now we can take $k'=\frac{k+\alpha-1+\beta}{2} \in \left( \beta,k+\alpha-1\right)$
so that
\begin{align}\label{21}
 1-\theta-\frac{ \theta \beta}{k+\alpha-1}=0.
\end{align}
From (\ref{inter}) and \er{xiaoyu}, using Young's inequality one obtains
\begin{align}
C(k)\|u\|_{L^{k'}(\R^n)}^{b} & \le C(k) \left( \|u\|_{L^{k+\alpha-1}(\R^n)}^{k+\alpha-1} \|u\|_{L^\beta(\R^n)}^\beta \right)^{\frac{b \theta}{k+\alpha-1}} \nonumber\\
& \le \frac{k}{4}  \|u\|_{L^{k+\alpha-1}(\R^n)}^{k+\alpha-1} \|u\|_{L^\beta(\R^n)}^\beta + C_0(k). \label{guji7}
\end{align}
Therefore, together with \er{nge3priori} and \er{n12priori} one has that for any $k>\max\left\{ \frac{2(\alpha-1)}{p-2},1 \right\}$ with $p$ is showed as in \er{p},
\begin{align}\label{globalesti}
\begin{array}{ll}
 \frac{d}{dt} \int_{\R^n}  u^k dx + \frac{3k}{4} \int_{\R^n} u^\beta dx \int_{\R^n}
u^{k+\alpha-1} dx + \frac{3(k-1)}{k} \left\| \nabla u^{\frac{k}{2}} \right\|_{L^{2}(\R^n)}^2
\le C_0(k),~~n \ge 3. \\
 \frac{d}{dt} \int_{\R^n}  u^k dx + \frac{3k}{4} \int_{\R^n} u^\beta dx \int_{\R^n}
u^{k+\alpha-1} dx + \frac{3(k-1)}{k} \left\| \nabla u^{\frac{k}{2}} \right\|_{L^{2}(\R^n)}^2 \\
\le C_0(k)+\frac{1}{2} \|u\|_{L^k(\R^n)}^k,~~n=1,2.
\end{array}
\end{align}
{\bf Step 2 ($L^k$ estimates for $\beta \le k \le \beta+\alpha-1$).} Based on the above a priori estimates, we firstly use the nonnegative term $\int_{\R^n} u^\beta dx \int_{\R^n} u^{\beta+\alpha-1} dx$ to get the uniform boundedness of $L^{\beta+\alpha-1}$ norm, as a consequence, it follows the estimates of $L^\beta$ norm.

First of all, by H\"{o}lder's inequality one has
\begin{align}
\|u\|_{L^{\beta+\alpha-1}(\R^n)}^{\beta+\alpha-1} & \le  \left( \|u\|_{L^{\beta+2(\alpha-1)}(\R^n)}^{\beta+2(\alpha-1)} \|u\|_{L^\beta(\R^n)}^{\beta} \right)^{\frac{1}{2}}.
\end{align}
Then by virtue of Young's inequality one has
\begin{align}
\|u\|_{L^{\beta+\alpha-1}(\R^n)}^{\beta+\alpha-1} \le \frac{\beta+\alpha-1}{2}\|u\|_{L^{\beta+2(\alpha-1)}(\R^n)}^{\beta+2(\alpha-1)} \|u\|_{L^\beta(\R^n)}^{\beta}+\frac{1}{2(\beta+\alpha-1)}. \label{kpline}
\end{align}
From \er{xiaoyup} and $p>2$ we know $\beta+\alpha-1>\max\left\{ \frac{2(\alpha-1)}{p-2},1 \right\}$, and then letting $k=\beta+\alpha-1$ in \er{globalesti} one has that for $n \ge 1$
\begin{align}
\frac{d}{dt} \int_{\R^n} u^{\beta+\alpha-1} dx + \int_{\R^n} u^{\alpha+\beta-1} dx \le C(\beta+\alpha-1),
\end{align}
which assures the following uniform estimate in time
\begin{align}\label{alphabeta}
\int_{\R^n} u^{\beta+\alpha-1} dx \le \|u_0\|_{L^{\alpha+\beta-1}(\R^n)}^{\alpha+\beta-1} e^{-t} +C(\alpha,\beta).
\end{align}

Next, we will apply the boundedness of $\int_{\R^n} u^{\beta+\alpha-1} dx$ norm to show that the $L^\beta$ norm is also uniformly bounded in time. By taking $k=\beta$ in \er{guji2} one obtains
\begin{align}\label{inequality}
\frac{d}{dt} \int_{\R^n} u^\beta dx \le \beta \int_{\R^n} u^{\beta+\alpha-1} dx \left( 1-\int_{\R^n} u^\beta dx \right).
\end{align}
If $\int_{\R^n} u_0^\beta dx \le 1$, then either $\int_{\R^n} u(t)^\beta dx \le 1$ for all $t>0$ or there exists a time interval $[t_0,t_0+\varepsilon)$ such that
$\int_{\R^n} u(t_0)^\beta dx=1$ and $\int_{\R^n} u(t)^\beta dx>1$ and increases for $t_0 \le t<t_0+\varepsilon$. On the other hand,
\begin{align*}
\frac{d}{dt} \int_{\R^n} u^\beta dx<0,~t_0<t<t_0+\varepsilon,
\end{align*}
which is a contradiction with the increasing of $\int_{\R^n} u(t)^\beta dx$ within $t_0 \le t<t_0+\varepsilon$. Therefore, $\int_{\R^n} u(t)^\beta dx \le 1$ for any $t \ge 0$.

For $\int_{\R^n} u_0^\beta dx > 1$, if $\int_{\R^n} u(t)^\beta dx>1$ for all $t>0$, then $\frac{d}{dt} \int_{\R^n} u^\beta dx<0$ and thus $\int_{\R^n} u^\beta dx<\int_{\R^n} u_0^\beta dx$. Otherwise, denote $t_0$ to be the first time such that $\int_{\R^n} u(t_0)^\beta dx=1$, then using the above arguments for $\int_{\R^n} u_0^\beta dx \le 1$ we know that
\begin{align*}
\frac{d}{dt} \int_{\R^n} u^\beta dx<0,~0 \le t<t_0, ~~\mbox{and}
\int_{\R^n} u^\beta dx \le 1,~t \ge t_0.
\end{align*}
Collecting the two cases we obtain
\begin{align}
\int_{\R^n} u^\beta dx \le \max \left\{ \int_{\R^n} u_0^\beta dx,1\right\}.
\end{align}
Hence, H\"{o}lder's inequality gives that for any $\beta \le k \le \beta+\alpha-1$, we have
\begin{align}
\int_{\R^n} u^k dx \le C(k,u_0).
\end{align}
{\bf Step 3 ($L^k$ estimates for $\beta+\alpha-1<k<\infty$).} In this step, based on \er{guji2}, the nonnegative term $\|\nabla u^{\frac{k}{2}}\|_{L^2(\R^n)}^2$ will be taken into account to obtain the boundedness of $L^k$ norm for $k>\beta+\alpha-1$. Letting
$$
v=u^{\frac{k}{2}},~q=2,~r=1<m<2,~C_0=\frac{k-1}{k},~C_1=\frac{1}{2}
$$
in Lemma \ref{interpolation} one has
\begin{align}\label{39}
\begin{array}{ll}
\int_{\R^n} u^{k} dx \le \frac{k-1}{k} \left\|  \nabla u^{\frac{k}{2}}
\right\|_{L^2(\R^n)}^{2 } + C(n,k) \|u\|_{L^{k_1}(\R^n)}^k,~~n \ge 3, \\
\int_{\R^n} u^{k} dx \le \frac{k-1}{k} \left\|  \nabla u^{\frac{k}{2}}
\right\|_{L^2(\R^n)}^{2 }+\frac{1}{2}\|u\|_{L^k(\R^n)}^k+ C(n,k) \|u\|_{L^{k_1}(\R^n)}^k,~~n=1,2, \\
\end{array}
\end{align}
where $k_1=\frac{km}{2}<k$. We can unify that for $n \ge 1$,
\begin{align}
\int_{\R^n} u^{k} dx \le \frac{2(k-1)}{k} \left\|  \nabla u^{\frac{k}{2}}
\right\|_{L^2(\R^n)}^{2 } + C(n,k) \|u\|_{L^{k_1}(\R^n)}^k.
\end{align}
On the other hand, by H\"{o}lder's inequality, for $k>\beta+\alpha-1$, we can take $k_1=\frac{\beta+k+\alpha-1}{2}<k$ so that
\begin{align}
\|u\|_{L^{k_1}(\R^n)}^k \le
\left( \|u\|_{L^{k+\alpha-1}(\R^n)}^{k+\alpha-1} \|u\|_{L^\beta(\R^n)}^\beta \right)^{\frac{k}{\beta+k+\alpha-1}}.
\end{align}
Hence it follows from Young's inequality that
\begin{align}
\frac{3}{2}\int_{\R^n} u^{k} dx & \le \frac{3(k-1)}{k} \left\|  \nabla u^{\frac{k}{2}}
\right\|_{L^2(\R^n)}^{2 } + C(n,k) \left( \|u\|_{L^{k+\alpha-1}(\R^n)}^{k+\alpha-1} \|u\|_{L^\beta(\R^n)}^\beta \right)^{\frac{k}{\beta+k+\alpha-1}} \nonumber \\
& \le  \frac{3(k-1)}{k} \left\|  \nabla u^{\frac{k}{2}}
\right\|_{L^2(\R^n)}^{2 } + \frac{3k}{4}  \|u\|_{L^{k+\alpha-1}(\R^n)}^{k+\alpha-1} \|u\|_{L^\beta(\R^n)}^\beta +C(n,k). \label{38}
\end{align}
Substituting \er{38} into \er{globalesti} we can have that for $n \ge 1$ and $k>\beta+\alpha-1$
\begin{align}
\frac{d}{dt} \int_{\R^n}  u^k dx + \int_{\R^n} u^{k} dx \le C(n,k).
\end{align}
Simple computations arrives that for any $k>\beta+\alpha-1$,
\begin{align}\label{44}
\int_{\R^n}  u^k dx \le \|u_0\|_{L^k(\R^n)}^k  e^{-t}+C(n,k).
\end{align}

{\bf Step 4 (The $L^\infty$ estimates).} On account of the above arguments, our last task is to give the uniform boundedness of solutions for any $t>0$.

Denote $q_k=2^k+\beta+\alpha-1$, by taking $k=q_k$ in \er{guji2}, we have
\begin{align}
\frac{d}{dt} \int_{\R^n}  u^{q_k} dx + \frac{4(q_k-1)}{q_k} \int_{\R^n}
|\nabla u^{\frac{q_k}{2}} |^2 dx+ q_k \int_{\R^n} u^\beta dx \int_{\R^n} u^{q_k+\alpha-1}
dx = q_k \int_{\R^n} u^{q_k+\alpha-1} dx. \label{45}
\end{align}
Armed with Lemma \ref{interpolation}, letting
$$
v=u^{\frac{q_k}{2}},~~q=\frac{2(q_k+\alpha-1)}{q_k},~~r=\frac{2 q_{k-1}}{q_k},~~C_0=\frac{1}{2 q_k},
$$
one has that for $n \ge 1$,
\begin{align}
\|u\|_{L^{q_k+\alpha-1}(\R^n)}^{q_k+\alpha-1} \le C(n) C_0^{\frac{-1}{\delta_1-1}} \left( \int_{\R^n} u^{q_{k-1}} dx \right)^{\gamma_1}+\frac{1}{2 q_k} \|\nabla u^{\frac{q_k}{2}} \|_{L^2(\R^n)}^2 + \frac{1}{2 q_k}\|u\|_{L^{q_k}(\R^n)}^{q_k},
\end{align}
where
\begin{eqnarray*}
\gamma_1&=&1+\frac{q_k+\alpha-1-q_{k-1}}{q_{k-1}-\frac{p(\alpha-1)}{p-2}} \le 2, \mbox{ iff } \alpha\leq 1+\Big(1-\frac{2}{p}\Big)\beta,\\
\delta_1&=&\frac{q_k- 2 q_{k-1}/p}{q_k+\alpha-1-q_{k-1}}=O(1).
\end{eqnarray*}
Substituting it into \er{45} and with notice that $\frac{4(q_k -1)}{q_k} \ge 2$, it follows
\begin{align}\label{47}
&\frac{d}{dt} \int_{\R^n} u^{q_k} dx +\frac{3}{2} \int_{\R^n} |\nabla u^{\frac{q_k}{2}}|^2 dx+ q_k \int_{\R^n} u^\beta dx \int_{\R^n} u^{q_k+\alpha-1}
dx \nonumber\\
\le & ~C(n) q_k^{\frac{\delta_1}{\delta_1 -1}} \left( \int_{\R^n} u^{q_{k-1}} dx \right)^{\gamma_1}+ \frac{1}{2} \|u\|_{L^{q_k}(\R^n)}^{q_k}.
\end{align}
Applying lemma \ref{interpolation} with
$$
v=u^{\frac{q_k}{2}},~~q=2,~~r=\frac{2 q_{k-1}}{q_k}
$$
and using Young's inequality, we have
\begin{align}\label{48}
& \frac{1}{2}\int_{\R^n} u^{q_k} dx \le C(n) \left( \int_{\R^n} u^{q_{k-1}} dx \right)^{\gamma_2} + \frac{1}{2}\int_{\R^n} |\nabla u^{\frac{q_k}{2}}|^2 dx \nonumber \\
\le &~\frac{\beta+\alpha-1}{3} \int_{\R^n} u^\beta dx \int_{\R^n}u^{q_k+\alpha-1}dx+C(n,\alpha,\beta)+ \frac{1}{2}\int_{\R^n} |\nabla u^{\frac{q_k}{2}}|^2 dx.
\end{align}
where we have used
$$
\gamma_2=1+\frac{q_k-q_{k-1}}{q_{k-1}} < 2, ~~q_{k-1}=\frac{q_k+\beta+\alpha-1}{2}.
$$
By summing up \er{47} and \er{48}, with the fact that $\gamma_1 \le 2$, we have
\begin{align*}
& \frac{d}{dt} \int_{\R^n} u^{q_k} dx+ \int_{\R^n} u^{q_k} dx  \le C(\delta_1) q_k^{\frac{\delta_1}{\delta_1 -1}} \left( \int_{\R^n} u^{q_{k-1}} dx \right)^{\gamma_1}+C(n,\alpha,\beta)\\
\le & \max [C(\delta_1),C(n,\alpha,\beta)] q_k^{\frac{\delta_1}{\delta_1 -1}} \left[ \left( \int_{\R^n} u^{q_{k-1}} dx \right)^{\gamma_1}+ 1 \right]\\
\le & 2 \max[C(\delta_1),C(n,\alpha,\beta)] q_k^{\frac{\delta_1}{\delta_1 -1}} \max\left\{ \left( \int_{\R^n} u^{q_{k-1}} dx \right)^{2},1\right\}.
\end{align*}

Let $K_0=\max\left\{1,\|u_0\|_{L^\beta(\R^n)},\|u_0\|_{L^\infty(\R^n)}\right\}$, we have the following inequality for initial data
\begin{align}
\int_{\R^n} u_0^{q_k} dx \le \left( \max \left\{ \|u_0\|_{L^\beta(\R^n)},\|u_0\|_{L^\infty(\R^n)} \right\} \right)^{q_k} \le K_0^{q_k}.
\end{align}
Let $d_0=\frac{\delta_1}{\delta_1 -1}$, it is easy to see that $q_k^{d_0}=[2^k+\beta+\alpha-1]^{d_0} \le [2^k+2^k(\beta+\alpha-1)]^{d_0}$. By taking $\bar{a}=2 \max \{C(\delta_1),C(n,\alpha,\beta)\} (\beta+\alpha)^{d_0}$ in the lemma 4.1 of \cite{BL14}, we obtain
\begin{align}\label{50}
\int_{\R^n} u^{q_k} dx \le (2 \bar{a})^{2^k-1} 2^{d_0 (2^{k+1}-k-2)} \max \left\{ \sup_{t\ge 0} \left(\int_{\R^n} u(t)^{q_0} dx\right)^{2^k}, K_0^{q_k} \right\}.
\end{align}

Since $q_k=2^k+\beta+\alpha-1$ and taking the power $\frac{1}{q_k}$ to both sides of \er{50}, then the boundedness of the solution $u(x,t)$ is obtained by passing to the limit $k \to \infty$
\begin{align}
\|u(t)\|_{L^\infty(\R^n)} \le 2 \bar{a} 2^{2 d_0} \max \left\{ \sup_{t\ge 0} \int_{\R^n} u(t)^{q_0} dx, K_0 \right\}.
\end{align}
On the other hand, by \er{44} with $q_0>\beta+\alpha-1$, we know
$$
\int_{\R^n} u(t)^{q_0} dx=\int_{\R^n} u(t)^{\beta+\alpha} dx \le \|u_0\|_{L^{\beta+\alpha}(\R^n)}^{\beta+\alpha}+C(\alpha,\beta) \le K_0^{\beta+\alpha}+C(\alpha,\beta).
$$
Therefore we finally have
\begin{align}
\|u(t)\|_{L^\infty(\R^n)} \le C(\alpha,\beta,K_0).
\end{align}
thus closes the proof. $\Box$

We now have necessary {\it a priori} estimates for the existence of global classical solutions and we know that $u$ is uniformly bounded for any $t \ge 0$. Moreover, the reaction term $u^\alpha \left(1-\int_{\R^n} u^\beta dx\right)$ is bounded from below and above. Hence the uniqueness and global existence of classical solution is followed by the standard parabolic theory. This completes the proof of Theorem \ref{thm1}.

\section{Conclusions}
This paper proves the global existence and uniqueness of solution to Cauchy problem \er{nkpp}. Especially when $\beta=1$, for $1 \le \alpha< 1+2/n$ with $n \ge 2$, it has been proved that $u^\alpha (1-\sigma\int_{\R^n} u dx)$ is bounded from below and above, therefore, if $u^\alpha (1-\sigma\int_{\R^n} u dx)$ is positive, the structure is similar to Fujita equation $u_t=\Delta u+u^\alpha$ in the whole space. However, our equation has global solution for $1 < \alpha<1+2/n$ when $ n\ge 2$, which is opposite to the result of Fujita equation, where no global solution exists for any initial data. The difference arises from the nonlocal term $1-\sigma\int_{\R^n} u dx$. In the modeling of population dynamics, there are more and more nonlocal reaction diffusion equations which have been derived. However, the corresponding mathematical theory is far from complete. This paper gives a first step analysis in studying what  have the nonlocal effects brought into the problem.



\end{document}